\documentstyle[11pt,amssymb,amstex]{amsart}

\newtheorem{thm}{Theorem}[section]
\newtheorem{lem}[thm]{Lemma}
\newtheorem{cor}[thm]{Corollary}

\newtheorem{defin}[thm]{Definition}
\newtheorem{rem}[thm]{Remark}
\newtheorem{prob}[thm]{Problem}

\def\a{\alpha}

  \def\G{\Gamma}

\begin{document}

\title[On homotopy of volterrian q.s.o.]
{On homotopy of volterrian quadratic stochastic operators}
\author{Farrukh Mukhamedov}
\address{Farrukh Mukhamedov\\
 Department of Computational \& Theoretical Sciences\\
Faculty of Sciences, International Islamic University Malaysia\\
P.O. Box, 141, 25710, Kuantan\\
Pahang, Malaysia} \email{{\tt far75m@@yandex.ru}}

\author{Mansoor Saburov}
\address{Mansoor Saburov\\
Department of Mechanics \& Mathematics\\
National University of Uzbekistan\\
Vuzgorodok, 100174 Tashkent, Uzbekistan} \email{{\tt
msaburov@@gmail.com}}

\begin{abstract}
In the present paper we introduce a notion of homotopy of two
Volterra operators which is related to fixed points of such
operators. It is establish a criterion when two Volterra operators
are homotopic, as a consequence we obtain that the corresponding
tournaments of that operators are the same. This, due to
\cite{Ga1}, gives us a possibility to know some information about
the trajectory of homotopic Volterra operators. Moreover, it is
shown that any Volterra q.s.o. given on a face has at least two
homotopic extension to the whole simplex.
\end{abstract}

\maketitle

\footnotetext[1]{The first author (F.M.) is on leave from National
University of Uzbekistan, Vuzgorodok, 100174 Tashkent, Uzbekistan}

\section{Introduction}

Since Lotka and Volterra's seminal and pioneering
works\cite{L,V1,V2} many decades ago, modeling of interacting,
competing species have received considerable attention in the
fields of biology, ecology, mathematics \cite{H,HS,RMH,T,MGT} and,
more recently, in the physics literature as well
\cite{DP,KLF,Lip1,Lip2,NSE,PL}. In their remarkably simple
deterministic model, Lotka and Volterra considered two coupled
nonlinear differential equations that mimic the temporal evolution
of a two-species system of competing predator and prey
populations. They demonstrated that coexistence of both species
was not only possible but inevitable in their model. Moreover,
similar to observations in real populations, both predator and
prey densities in this deterministic system display regular
oscillations in time, with both the amplitude and the period
determined by the prescribed initial conditions.

While investigation of computational aspects of such dynamical
systems, it needs to consider discretization of such systems. This
leads to study of the trajectory of discrete time Volterra
operators. Therefore, in \cite{Ga1,Ga2,GE,M,MAT}\footnote[2]{Note
that more general, so {\it called quadratic operators}, were
studied by many authors (see for example,\cite{B,K, Lyu})}
discrete time Volterra operators were considered and investigated.
A connection between such dynamical systems and the theory of
tournaments were established. This gave some information about the
trajectory of Volterra operators, since the corresponding
tournaments are related to the fixed points of Volterra operators.
Moreover, some ergodic properties of such operators, in small
dimensions, were studied in \cite{GZ,Val,Z}. However, still much
information is unknown about behavior of Volterra operators.

In the present paper we introduce a notion of homotopy of two
Volterra operators which is related to fixed points of  such
operators. Further, we will establish a criterion when two
Volterra operators are homotopic, as a consequence we obtain that
the corresponding tournaments of that operators are the same.
This, due to \cite{Ga1}, gives us a possibility to know some
information about the trajectory of homotopic Volterra operators.
Moreover, it is shown that any Volterra q.s.o. given on a face has
at least two homotopic extension to the whole simplex.

\section{Preliminaries}

Denote by $$S^{m-1}=\left\{x=(x_1,x_2,...,x_m)\in R^m:
\sum\limits_{k=1}^mx_k=1, \ x_k\ge 0)\right\}$$ $(m-1)-$
dimensional simplex. The vertices of the simplex $S^{m-1}$ are
described by the elements
$e_k=(\delta_{1k},\delta_{2k},\dots,\delta_{mk})$, where
$\delta_{ik}-$ is the Kronecker's symbol. Let $I=\{1,2,\dots m\}$
and $\alpha\subset I-$ be an arbitrary subset. By $\Gamma_\alpha$
we denote the convex hull of the vertices $\{e_i\}_{i\in\alpha}.$
The set $\Gamma_\alpha$ is usually called
$(|\alpha|-1)-$dimensional face of the simplex, here $|\alpha|$
stands for the cardinality of  $\alpha$. An interior of
$\Gamma_\alpha$ in the induced topology of ${\mathbb{R}}^m$ to
affine hull $\Gamma_\alpha$ is called relative interior and is
denoted by $ri\Gamma_\alpha.$ One can see that that
$$ri\Gamma_\alpha=\Bigl\{x\in S^{m-1}:  x_k>0 \ \forall k\in\alpha; \ x_k=0 \ \forall k\notin\alpha \Bigr\}.$$
Similarly, one can define relative boundary
$\partial\Gamma_\alpha$ of the face $\Gamma_\alpha.$ In
particular, we have
$$riS^{m-1}=\{x\in S^{m-1}: x_k>0 \ \forall k\in I\},$$
$$\partial S^{m-1}=\{x\in S^{m-1}: \exists k\in I; x_k=0\}.$$

A {\it Volterra quadratic stochastic operator (q.s.o.)}
$V:S^{m-1}\rightarrow S^{m-1}$ is defined by
\begin{equation}\label{1}
(V(x))_k=x_k\left(1+\sum\limits_{i=1}^{m}a_{ki}x_i\right), \
k=\overline{1,m},
\end{equation}
where $a_{ki}=-a_{ik}$, $|a_{ki}|\le 1,$ i.e.
$A_m=\left(a_{ki}\right)_{k,i=1}^m$ is a skew-symmetrical matrix.

Let
$$Fix(V)=\{x\in S^{m-1}: Vx=x\}$$
be the set of all fixed points of the Volterra q.s.o. $V$.

One can see that for any Volterra operator $V$ the set $Fix(V)$
contains all the vertices of the simplex $S^{m-1}$. Therefore, it
is non empty.

For given $x\in S^{m-1}$ consider a sequence $\{x,Vx,\dots,
V^nx,\dots\}$ which is called {\it the trajectory} of a Volterra
q.s.o. $V$. Limit points of such a sequence is denoted by
$\omega_V(x)$.

 For an arbitrary $x\in
{\mathbb{R}}^m$ let us define
$$supp(x)=\{i\in I: x_i\neq 0\}$$
\textit{support} of the $x$. The following statements can be
easily proved for Volterra q.s.o.

\begin{thm}\label{vol1} \cite{Ga1} Let $V:S^{m-1}\rightarrow S^{m-1}$
be a Volterra q.s.o. and $V_\alpha-$ be the restriction $V$ to the
face $\Gamma_\alpha$, then the following assertions hold true:
\begin{itemize}
    \item[(i)] For any $\alpha\subset I$ one has
    $V(\Gamma_\alpha)\subset\Gamma_\alpha$.

    \item[(ii)] For any $\alpha\subset I$ one has
    $V(ri\Gamma_\alpha)\subset ri\Gamma_\alpha$ and $V(\partial\Gamma_\alpha)\subset \partial\Gamma_\alpha.$

    \item[(iii)] The restriction $V_\alpha:\Gamma_\alpha\rightarrow\Gamma_\alpha$ is also a Volterra q.s.o.

    \item[(iv)] If $x\in Fix(V)X$ then $Supp(x) \cap Supp(A_mx)=\emptyset,$ in particular if
    $x\in Fix(V_\cap riS^{m-1}$ then $x\in Ker A_m,$, here $Ker A_m$ is the kernel of the matrix $A_m.$

    \item[(v)] For any $x\in S^{m-1}$ the set $\omega_V(x)$ either consists of a single point or is infinite.

     \item[(vi)] The set of all volterra q.s.o.s geometrically can be considered
    as a      $\left(\frac{m(m-1)}{2}\right)-$ dimensional cub on ${\mathbb{R}}^{\frac{m(m-1)}{2}}.$
\end{itemize}
\end{thm}

Let  $A_m$ be a skew-symmetrical matrix corresponding to a
Volterra q.s.o. given by \eqref{1}. It is known \cite{PR} that if
the order of a skew-symmetrical matrix is odd, then the
determinant of this matrix is 0, otherwise the determinant is the
square of some polynomial of its entries, which situated above the
main diagonal. Such a polynomial is called \textit{pffaffian} and
can by calculated by the following rule.

\begin{lem}[\cite{PR}]\label{pfaf} Let $p_m-$ be a pffaffian of an even order skew-symmerical matrix $A_m$
$m>2$. By $p_{im}-$ we denote a pffaffian of the skew-symmerical
matrix $A_{im}$, which is obtained from $A_m$ by deleting of the
$m$-th and $i$-th rows and columns, where $i=\overline{1,m-1}$.
Then one has
$$p_m=\sum\limits_{i=1}^{m-1}(-1)^{i-1}p_{im}a_{im}; \ \ \ p_2=a_{12},
$$
where $p_{im}$ is obtained  via $p_{m-2}$ by adding 1 to all
indexes greater or equal to $i$.
\end{lem}

A pffaffian of the main minor of an even order skew-symmetric
matrix $A_m$ with rows and columns $\{i_1,i_2,\dots,i_{2k}\}$ is
called {\it main subpffafian} of order $2k$, and is denoted by
$gp_{i_1i_2\dots i_{2k}}.$

For example,  $$gp_{i_1i_2}=a_{i_1i_2};  \ \
gp_{i_1i_2i_3i_4}=a_{i_1i_2}a_{i_3i_4}+a_{i_1i_4}a_{i_2i_3}-a_{i_1i_3}a_{i_2i_4}.$$

\begin{defin}\label{11} A skew-symmetrical matrix $A_m$ is called
transversal if all even order main minors are nonzero.
\end{defin}\

It is clear that if a skew-symmetrical matrix $A_m$ is
transversal, then all even order main subpffaffians are nonzero,
that is
$$gp_{i_1i_2\dots i_{2k}}\neq 0, \ \ \forall i_1,i_2,\dots,i_{2k}\in I,$$
in particular $a_{ki}\neq 0$ for $k\neq i.$

\begin{defin}\label{2} We say that a Volterra q.s.o. $V$ is transversal if the
corresponding skew-symmetrical matrix $A_m$ is transversal.
\end{defin}

Denote by ${\mathcal{V}}_t^{m-1}$ the set of all transversal
Volterra q.s.o.s defined in the simplex $S^{m-1}.$

Henceforth, we will consider only transversal Volterra q.s.o., and
do not use the word "transversal".

\begin{thm}\label{vol2}\cite{Ga1} Let $V\in{\mathcal{V}}_t^{m-1},$ then
\begin{itemize}
    \item[(i)] $Fix(V)$ is a finite set,
    \item[(ii)] If $x\in Fix(V)$, then the cardinality of $Supp (x)$ is odd.
    \item[(iii)] For any face $\Gamma_\alpha$ of the simplex $S^{m-1}$ one has
    $|Fix(V)\cap ri\Gamma_\alpha|\le 1.$
\end{itemize}
\end{thm}

\begin{rem}\label{fix} Note that there is no fixed points
of any  (transversal) Volterra  q.s.o. in the interior of odd
dimensional faces.
\end{rem}

\section{Homotopy of Volterra operators}

In this section we are going to define a notion of homotopy for
Volterra q.s.o.s. Further, we will show that two homotopic
Volterra q.s.o. have 'similar' trajectories under some conditions.

\begin{defin}\label{3} Two Volterra q.s.o. $V_0, V_1\in{\mathcal{V}}_t^{m-1}$ are  called
homotopic, if there exists a family of Volterra operators
$\{V_\lambda\}_{\lambda\in [0,1]}\subset {\mathcal{V}}_t^{m-1}$
such that it is continuous with respect to $\lambda$ with
$V_\lambda\mid_{\lambda=0}=V_0$, $V_\lambda\mid_{\lambda=1}=V_1$
and one has $|Fix(V_\lambda)|=|Fix(V_0)|=|Fix(V_1)|$ for any
$\lambda\in [0,1]$.
\end{defin}

\begin{rem}\label{pfaf} Note that if a family $\{V_\lambda\}_{\lambda\in [0,1]}\subset {\mathcal{V}}_t^{m-1}$
is continuous than one can see that the main subpffaffians
$gp_{i_1i_2\dots i_{2k}}^{(\lambda)}$ of the corresponding
skew-symmetric matrices $A_m^{(\lambda)}$ are also continuous with
respect to $\lambda$.
\end{rem}

One can see that the introduced homotopy defines an equivalency
relation in the set ${\mathcal{V}}_t^{m-1}$. Therefore, two
operators $V_0, V_1\in{\mathcal{V}}_t^{m-1}$ are called
\textit{equivalent} and denoted by $V_0\sim V_1,$ if they are
homotopic. Hence, one can consider a factor set
${\mathcal{V}}_t^{m-1}\diagup_\sim .$

{\bf Example.} Let $m=2$. Then Voterra operators corresponding to
the following  matrices
$$A_a=\left(%
\begin{array}{cc}
  0 & a \\
  -a & 0 \\
\end{array}
\right), \ \ \ \ 0<a\le 1$$ are always homopotic.

Let $m=3$. Then Volterra operators corresponding to the following
matrices
$$A_{abc}=\left(%
\begin{array}{ccc}
  0 & a & b \\
  -a & 0 & c \\
  -b & -c & 0 \\
\end{array}%
\right),\ \ \ \  0<a,b,c\le 1$$ are always homotopic.

 Now we are interested when two Volterra q.s.o. are
equivalent.

\begin{thm}\label{eqv1} Let $V_0, V_1\in{\mathcal{V}}_t^{m-1}$ with $V_0 \sim V_1$ and
$A_m^{(0)}$,$A_m^{(1)}$ be their the corresponding skew-symmetric
matrices. Then all corresponding even order main subpffaffians the
matrices $A_m^{(0)}$ and $A_m^{(1)}$ have the same sign i.e.
$$Sign (gp_{i_1i_2\dots i_{2k}}^{(0)})=Sign(gp_{i_1i_2\dots
i_{2k}}^{(1)})\ \
 \ \ \forall i_1,i_2,\dots ,i_{2k}\in I.$$
\end{thm}

\begin{pf} Due to $V_0 \sim V_1$ there exists a continuous family
$\{V_\lambda\}_{\lambda\in[0,1]}\in{\mathcal{V}}_t^{m-1}$ such
that  $V_\lambda\mid_{\lambda=0}=V_0$ and
$V_\lambda\mid_{\lambda=1}=V_1$.  Let us consider a
skew-symmetrical matrix $A_m^{(\lambda)}$ corresponding to
$V_\lambda$. Then $A_m^{(\lambda)}\mid_{\lambda=0}=A_m^{(0)}$ and
$A_m^{(\lambda)}\mid_{\lambda=1}=A_m^{(1)}.$

Assume that the assertion of the theorem is not true, that is
there are $2k_0$ order main subpffaffians of the matrices
$A_m^{(0)}$ and $A_m^{(1)}$ such that
$$Sign(gp_{i_1i_2\dots i_{2k_0}}^{(0)})\neq Sign
(gp_{i_1i_2\dots i_{2k_0}}^{(1)})$$ which implies
$$gp_{i_1i_2\dots i_{2k_0}}^{(0)}gp_{i_1i_2\dots i_{2k_0}}^{(1)}<0.
$$
Continuity of $gp_{i_1i_2\dots i_{2k_0}}^{(\lambda)}$ with respect
to $\lambda$ (see Remark \ref{pfaf}) yields the existence of
$\lambda_0\in [0,1]$ such that $gp_{i_1i_2\dots
i_{2k_0}}^{(\lambda_0)}=0.$ But the last contradicts to
$V_{\lambda_0}\in{\mathcal{V}}_t^{m-1}.$
\end{pf}

Let us recall some definitions relating to tournaments associated
with a skew-symmetrical matrix $A_m=(a_{ki})_{k,i=1}^m$. Put
$$
Sign(A_m)=\left(Sign \ a_{ki}\right)_{k,i=1}^m.
$$

Define a {\it tournament} $T_{m}$, as a graph consisting of $m$
vertices labelled by $\{1,2,\dots,m\}$, corresponding to a
skew-symmetrical matrix $A_m$ by the following rule: there is an
arrow from $i$ to $k$ if $a_{ki}<0$, a reverse arrow otherwise.
Note that if signs of two skew-symmetric matrices are the same,
then the corresponding tournaments are the same as well.

Recall that a tournament is said to be \textit{strong} if it is
possible to go from any vertex to any other vertex with directions
taken into account. A \textit{strong component} of a tournament is
a maximal strong subtournament of the tournament. The tournament
with the strong components of $T_m$ as vertices and with the edge
directions induced from $T_m$ is called \textit{the factor
tournament} of the tournament $T_m$ and denoted by
$\widetilde{T}_m.$ \textit{Transitivity} of the tournament means
that there is no strong subtournament consisting of three vertices
of the given tournament. A tournament containing fewer than three
vertices is regarded as \textit{transitive} by definition. As is
known \cite{HF}, the factor tournament $\widetilde{T}_m$ of any
tournament $T_m$ is transitive. Further, after a suitable
renumbering of the vertices of $T_m$ we can assume that the
subtournament $T_r$ contains the vertices of $T_m$ as its
vertices, i.e., $\{1\}, \{2\}, \cdots ,\{r\}.$ Obviously, $r\ge
m,$ and $r=m$ if and only if $T_m$ is a strong tournament.

\begin{cor}\label{Tour1} If $V_0 \sim V_1,$ then the corresponding tournaments $T_m^{(0)}$
and $T_m^{(1)}$ are the same.
\end{cor}

\begin{pf}
Since, $gp_{ki}=a_{ki}$ then Theorem \ref{eqv1} implies that $Sign
(A_m^{(0)})=Sign(A_m^{(1)}).$ Hence, the corresponding tournaments
$T_m^{(0)},$ $T_m^{(1)}$ are the same.
\end{pf}

This Corollary gives some information about the trajectory of
equivalent Volttera operators. Namely, due to results of
\cite{Ga1} and Corollary \ref{Tour1} one gets the following

\begin{cor}\label{Tour2} Let $V_0 \sim V_1$. The following assertions hold true:
\begin{enumerate}
\item[(i)] Assume that the tournament $T_m^{(0)}$ corresponding to
$V_0$ is not strong. Then  for any $x^0\in intS^{m-1}$ and $i>r,$
then $\omega_{V_0}(x^0)\subset \Gamma_\alpha$, and
$\omega_{V_1}(x^0)\subset \Gamma_\alpha$, here
$\alpha=\{1,2,\dots, r\}.$ \item[(ii)] Assume that $T_m^{(0)}$ is
transitive, then for any $x^0\in riS^{n-1}$
$\omega_{V_0}(x^0)=\omega_{V_1}(x^0)=(1,0,\dots,0).$
\end{enumerate}
\end{cor}

In Theorem \ref{eqv1} we have formulated a necessary condition to
be equivalent of two Volterra q.s.o. Now in small dimensions, we
are going to provide certain criterions for equivalence.

\begin{thm}\label{eqv2} Let $m\le 3$. Then  $V_0 \sim V_1$ if and only if
$$Sign(A_m^{(0)})=Sign (A_m^{(1)}).$$
\end{thm}

\begin{pf}
Necessity immediately follows from Theorem \ref{eqv1}. Therefore,
we will prove sufficiency. Consider separately two distinct case
with respect to $m$.

Let $m=2.$ Then
$$ A_2^{(0)}=\left(%
\begin{array}{cc}
  0 & a_{12}^{(0)} \\
  -a_{12}^{(0)} & 0 \\
\end{array}%
\right) \ \ \mbox{è} \ \ A_2^{(1)}=\left(%
\begin{array}{cc}
  0 & a_{12}^{(1)} \\
  -a_{12}^{(1)} & 0 \\
\end{array}%
\right) $$ here $Sign \ a_{12}^{(0)}=Sign \ a_{12}^{(1)}.$

Consider $A_2^{(\lambda)}=(1-\lambda)A_2^{(0)}+\lambda A_2^{(1)}.$
It is clear that $A_2^{(\lambda)}$ is transversal for any
$\lambda\in [0,1].$ Let $V_\lambda$ be the corresponding Volterra
q.s.o. (see \eqref{1}). Then one has
$V_\lambda=(1-\lambda)V_0+\lambda V_1$ and
$\{V_\lambda\}_{\lambda\in[0,1]}\subset {\mathcal{V}}_t^1$.
According to Theorem \ref{vol2} and Remark \ref{fix} the set of
all fixed points of $V_\lambda\in{\mathcal{V}}_t^1$ consists of
the vertices of $S^1.$ Therefore,
$|Fix(V_\lambda)|=|Fix(V_0)|=|Fix(V_1)|=2$ for any $\lambda\in
[0,1].$

Let $m=3$. Then
$$A_3^{(0)}=\left(%
\begin{array}{ccc}
  0 & a_{12}^{(0)} & a_{13}^{(0)} \\
  -a_{12}^{(0)} & 0 & a_{23}^{(0)} \\
  -a_{13}^{(0)} & -a_{23}^{(0)} & 0 \\
\end{array}%
\right), \ \  \ \ A_3^{(1)}=\left(%
\begin{array}{ccc}
  0 & a_{12}^{(1)} & a_{13}^{(1)} \\
  -a_{12}^{(1)} & 0 & a_{23}^{(1)} \\
  -a_{13}^{(1)} & -a_{23}^{(1)} & 0 \\
\end{array}%
\right)$$ here $Sign \ a_{ij}^{(0)}=Sign \ a_{ij}^{(1)}$ for
$i<j.$

One can check that the skew-symmetrical matrix
$A_3^{(\lambda)}=(1-\lambda)A_3^{(0)}+\lambda A_3^{(1)}$ is
transversal for any $\lambda\in [0,1].$ Therefore, the
corresponding Volterra q.s.o. $V_\lambda$ belongs to
${\mathcal{V}}_t^2$ and one has  $V_\lambda=(1-\lambda)V_0+\lambda
V_1$ for any $\lambda\in [0,1].$

Assume that $V\in{\mathcal{V}}_t^2$,  and its corresponding matrix
be
$$A_3=\left(%
\begin{array}{ccc}
  0 & a_{12} & a_{13} \\
  -a_{12} & 0 & a_{23} \\
  -a_{13} & -a_{23} & 0 \\
\end{array}%
\right).$$ Then one can find that if
$$Sign (a_{12})=Sign
(a_{13})=Sign(a_{23}),
$$
then $|Fix(V)|=4$, otherwise $|Fix(V)|=3.$

Hence, if the condition of the theorem is satisfied i.e.
$Sign(A_3^{(0)})=Sign(A_3^{(1)})$ then from the last we conclude
that either $|Fix(V_0)|=|Fix(V_1)|=4$ or
$|Fix(V_0)|=|Fix(V_1)|=3.$

Due to  $Sign(A_3^{(\lambda)})=Sign(A_3^{(0)})=Sign(A_3^{(1)}),$
one gets $$|Fix(V_\lambda)|=|Fix(V_0)|=|Fix(V_1)|$$ for any
$\lambda\in [0,1].$
\end{pf}

\begin{cor}\label{factor} If $m=2,$ then $|{\mathcal{V}}_t^2\diagup_\sim|=2$ and
if $m=3,$ then $|{\mathcal{V}}_t^3\diagup_\sim|=8$
\end{cor}

\begin{thm}\label{eqv3} Let $m=4$. Then $V_0 \sim V_1$ if and only
if
$$Sign(A_4^{(0)})=Sign(A_4^{(1)}), \ \  Sign (gp_{1234}^{(0)})=Sign(gp_{1234}^{(1)}),$$
where
$$
gp_{1234}^{(i)}=a_{12}^{(i)}a_{34}^{(i)}+a_{14}^{(i)}a_{23}^{(i)}-a_{13}^{(i)}a_{24}^{(i)},
\ \ i=0,1.
$$
\end{thm}

\begin{pf} As before, the necessity immediately follows from Theorem \ref{eqv1}.
Let us prove the sufficiency.

Let
$$A_4^{(0)}=\left(%
\begin{array}{cccc}
  0 & a_{12}^{(0)} & a_{13}^{(0)} & a_{14}^{(0)} \\
  -a_{12}^{(0)} & 0 & a_{23}^{(0)} & a_{24}^{(0)} \\
  -a_{13}^{(0)} & -a_{23}^{(0)} & 0 & a_{34}^{(0)} \\
  -a_{14}^{(0)} & -a_{24}^{(0)} & -a_{34}^{(0)} & 0 \\
\end{array}%
\right), \ \
A_4^{(1)}=\left(%
\begin{array}{cccc}
  0 & a_{12}^{(1)} & a_{13}^{(1)} & a_{14}^{(1)} \\
  -a_{12}^{(1)} & 0 & a_{23}^{(1)} & a_{24}^{(1)} \\
  -a_{13}^{(1)} & -a_{23}^{(1)} & 0 & a_{34}^{(1)} \\
  -a_{14}^{(1)} & -a_{24}^{(1)} & -a_{34}^{(1)} & 0 \\
\end{array}%
\right);
$$ here $Sign \ a_{ij}^{(0)}=Sign \ a_{ij}^{(1)}$ for $i<j$ and
$$Sign\left(a_{12}^{(0)}a_{34}^{(0)}+a_{14}^{(0)}a_{23}^{(0)}-a_{13}^{(0)}a_{24}^{(0)}\right)
=Sign\left(a_{12}^{(1)}a_{34}^{(1)}+a_{14}^{(1)}a_{23}^{(1)}-a_{13}^{(1)}a_{24}^{(1)}\right).$$

Let us consider the following skew-symmetric matrix
$A_4^{(\lambda)}$ defined by
$$\left(%
\begin{array}{cccc}
  0 & (1-\lambda)a_{12}^{(0)}+\lambda a_{12}^{(1)} & (1-\lambda)a_{13}^{(0)}+\lambda a_{13}^{(1)} & (1-\lambda)a_{14}^{(0)}+\lambda a_{14}^{(1)} \\
  -\left((1-\lambda)a_{12}^{(0)}+\lambda a_{12}^{(1)}\right) & 0 & \frac{(1-\lambda)a_{14}^{(0)}a_{23}^{(0)}+\lambda a_{14}^{(1)}a_{23}^{(1)}}{(1-\lambda)a_{14}^{(0)}+\lambda a_{14}^{(1)}} & \frac{(1-\lambda)a_{13}^{(0)}a_{24}^{(0)}+\lambda a_{13}^{(1)}a_{24}^{(1)}}{(1-\lambda)a_{13}^{(0)}+\lambda a_{13}^{(1)}} \\
  -\left((1-\lambda)a_{13}^{(0)}+\lambda a_{13}^{(1)}\right) & -\frac{(1-\lambda)a_{14}^{(0)}a_{23}^{(0)}+\lambda a_{14}^{(1)}a_{23}^{(1)}}{(1-\lambda)a_{14}^{(0)}+\lambda a_{14}^{(1)}} & 0 & \frac{(1-\lambda)a_{12}^{(0)}a_{34}^{(0)}+\lambda a_{12}^{(1)}a_{34}^{(1)}}{(1-\lambda)a_{12}^{(0)}+\lambda a_{12}^{(1)}} \\
  -\left((1-\lambda)a_{14}^{(0)}+\lambda a_{14}^{(1)}\right) & -\frac{(1-\lambda)a_{13}^{(0)}a_{24}^{(0)}+\lambda a_{13}^{(1)}a_{24}^{(1)}}{(1-\lambda)a_{13}^{(0)}+\lambda a_{13}^{(1)}} & -\frac{(1-\lambda)a_{12}^{(0)}a_{34}^{(0)}+\lambda a_{12}^{(1)}a_{34}^{(1)}}{(1-\lambda)a_{12}^{(0)}+\lambda a_{12}^{(1)}} & 0 \\
\end{array}%
\right).$$ It is then clear that  $Sign
(a_{ij}^{(\lambda)})=Sign(a_{ij}^{(0)})=Sign (a_{ij}^{(1)})$ for
$i<j$ and
\begin{eqnarray*}
Sign\bigg(a_{12}^{(\lambda)}a_{34}^{(\lambda)}+a_{14}^{(\lambda)}a_{23}^{(\lambda)}-a_{13}^{(\lambda)}a_{24}^{(\lambda)}\bigg)&=&
Sign\bigg((1-\lambda)a_{12}^{(0)}a_{34}^{(0)}+\lambda a_{12}^{(1)}a_{34}^{(1)}\\
&&+(1-\lambda)a_{14}^{(0)}a_{23}^{(0)}+\lambda a_{14}^{(1)}a_{23}^{(1)}\\
&&-(1-\lambda)a_{13}^{(0)}a_{24}^{(0)}-\lambda a_{13}^{(1)}a_{24}^{(1)}\bigg) \\
&=&Sign\left(a_{12}^{(0)}a_{34}^{(0)}+a_{14}^{(0)}a_{23}^{(0)}-a_{13}^{(0)}a_{24}^{(0)}\right)\\
&=&Sign\left(a_{12}^{(1)}a_{34}^{(1)}+a_{14}^{(1)}a_{23}^{(1)}-a_{13}^{(1)}a_{24}^{(1)}\right).
\end{eqnarray*}
for any $\lambda\in[0,1]$. This implies that the corresponding
Volterrian operator $V_\lambda$ belongs to ${\mathcal{V}}_t^4$ for
any $\lambda\in[0,1].$

Since $m=4$ is even, then according to Remark \ref{fix} there is
no fixed point in the interior of the simplex $S^3$. Thanks to
Theorem \ref{eqv2} one gets
$$
|Fix(V_\lambda) \cap
\partial S^3|=|Fix(V_0) \cap \partial S^3|=|Fix(V_1) \cap \partial S^3|
$$
for any $\lambda\in[0,1].$ Therefore,
$|Fix(V_\lambda)|=|Fix(V_0)|=|Fix(V_1)|$ for any
$\lambda\in[0,1].$
\end{pf}

\begin{cor}\label{factor2} If $m=4,$ then $|{\mathcal{V}}_t^4\diagup_\sim|= 112.$
\end{cor}

The following theorem can be considered a reverse to Theorem
\ref{eqv1}.

\begin{thm}\label{eqv4} Let $V_0, V_1\in {\mathcal{V}}_t^{m-1}$ and
$A_m^{(0)}$, $A_m^{(1)}$ be their the corresponding skew-symmetric
matrices. If all corresponding even order main subpffaffians of
the matrices
 $A_m^{(0)}$ and $A_m^{(1)}$ have the same sign, that is
$$Sign(gp_{i_1i_2\dots i_{2k}}^{(0)})=Sign(gp_{i_1i_2\dots i_{2k}}^{(1)}),
 \ \ \forall i_1,i_2,\dots,i_{2k}\in I,$$
then $|Fix(V_0)|=|Fix(V_1)|.$
\end{thm}

\begin{pf}
We will prove this by the induction with respect to the dimension
$m$ of the simplex $S^{m-1}$. For small dimensions our assumption
is true (see Theorems \ref{eqv2} and \ref{eqv3}). Let us assume
that the statement of the theorem is true for dimension $m-1$. Now
we prove it for dimension $m.$

Since the restriction of any transversal Volterra q.s.o. to any
face $\Gamma_\alpha$ of the simplex is also transversal Volterra
q.s.o. (see Theorem \ref{vol1}), then by the assumption of the
induction we get that operators $V_0$ and $V_1$ have the same
number of fixed points in $\partial S^{m-1}$, i.e.
$$Fix(V_0)\cap\partial S^{m-1}=Fix(V_1)\cap\partial S^{m-1}.$$

Let us show that the operators $V_0$ and $V_1$ have the same
number of fixed points in $ri S^{m-1}$.

If $m$ is even, then due to Remark \ref{fix} there is no any fixed
point of Volterra q.s.o. in the interior of the simplex.
Therefore, we have to prove the theorem only when $m$ is odd. Then
in this case, Theorem \ref{vol2} implies that $|Fix(V)\cap ri
S^{m-1}|\le 1$ for any $V\in{\mathcal{V}}_t^{m-1}$. According to
Theorem \ref{vol1} one can see that $x\in Fix(V)\cap ri S^{m-1}$
if and only if $x\in KerA_m\cap ri S^{m-1}$.

Due to oddness of $m$ the determinant of $A_m$ equals to 0, but
the transversality of the operator $V_0$ implies that the minor of
order $m-1$ is not zero, which means dimension of the image
$Im(A_m)$ is $m-1$. Hence, the equality $dim(KerA_m)+dim(ImA_m)=m$
implies that $Ker A_m$ is a one dimensional space.

Now we are going to describe $Ker A_m$. Keeping in mind that $det
A_m$ is zero, one finds
\begin{equation}\label{2}
\sum\limits_{k=1}^ma_{ki}A_{ki}=det A_m=0, \ \ \ \forall
k=\overline{1,m},
\end{equation}
here $A_{ki}$ is an algebraic completion (i.e. algebraic minor) of
entry $a_{ki}.$ It is known  \cite{PR} that
\begin{equation}\label{3}
A_{ki}=(-1)^{k+i}gp_{I_k}gp_{I_i}
\end{equation}
here as before $gp_{I_k}, gp_{I_i}$ are pffaffians of the minors
$M_{kk},M_{ii}$, where $I_k=I\setminus\{k\}, I_i=I\setminus\{i\}.$
It then follows from \eqref{2}, \eqref{3} that
$$
(-1)^kgp_{I_k}\sum\limits_{i=1}^ma_{ki}(-1)^igp_{I_i}=0, \ \ \
\forall k=\overline{1,m}.
$$ Thanks to $gp_{I_k}\neq 0,$ one finds
\begin{equation}\label{4}
\sum\limits_{i=1}^ma_{ki}(-1)^igp_{I_i}=0, \ \ \ \forall
k=\overline{1,m}.
\end{equation}
This means that for an element defined by
$$
x_0=\left(-gp_{I_1},gp_{I_2},\dots,(-1)^igp_{I_i},\dots
(-1)^mgp_{I_m}\right)
$$
one has $A_m(x_0)=0$, i.e. $x_0\in Ker A_m$. The
one-dimensionality of $A_m$ implies that $Ker A_m =\{\lambda x_0:
\ \lambda\in \mathbb{R}\}$. This means that there is an interior
fixed point for the Volterra operator $V_0$ if and only if
$$
Sign(-1)^1gp_{I_1}=Sign(-1)^2gp_{I_2}=\dots=Sign(-1)^ip_{I_i}=\dots=Sign(-1)^mgp_{I_m}$$
and that fixed point is given by
\begin{equation}\label{5}
x=\frac{1}{GP} x_0,
\end{equation} where
$$GP=\sum_{i=1}^{m}(-1)^kgp_{I_i}.$$

Now if the condition of the theorem is satisfied, then from
\eqref{5} one concludes that
$$|Fix(V_0)\cap ri S^{m-1}|=|Fix(V_1)\cap ri S^{m-1}|.
$$ Consequently, one gets $|Fix(V_0)|=|Fix(V_0)|$. This completes the proof.
\end{pf}

According to Theorem \ref{vol1} the set of all Volterra q.s.o.
geometrically forms a $\left(\frac{m(m-1)}{2}\right)$-dimensional
cube ${\mathcal{V}}^{m-1}$ in ${\mathbb{R}}^{\frac{m(m-1)}{2}}$.
Now let us consider the following manifolds
$$
\bigg\{V\in {\mathcal{V}}^{m-1}: gp_{i_1i_2\dots i_{2k}}(V)=0, \ \
\exists i_1,i_2,\dots,i_{2k}\in I\bigg\}.$$ These manifolds divide
the cube into several connected components.

From Theorems \ref{eqv1} and \ref{eqv4} one can prove the
following

\begin{thm}\label{eqv5} Two Volterra q.s.o. $V_0$ and $V_1$ ($V_0,V_1\in  {\mathcal{V}}_t^{m-1}$)
are homotopic if and only if the operators $V_0$ and $V_1$ belong
to only one connected component of the cube.
\end{thm}

\begin{pf} 'If' part of the proof immediately follows from Theorem
\ref{eqv1}. Therefore, let us prove 'only if' part.

Let us assume that $V_0$ and $V_1$ belong to the same connected
component. Then from the definition of component one can conclude
that such operators can be connected by a continuous path
$\{V_\lambda\}\in {\mathcal{V}}_t^{m-1}$ located in that
component.  On the other hand, we see that  the corresponding all
main subpffaffians of all operators $V_\lambda$ have the same
signs. So, thanks Theorem \ref{eqv4} one finds that
$|Fix(V_\lambda)|=|Fix(V_0)|=|Fix(V_1)|$ which implies that
$V_0\sim V_1$.
\end{pf}

\begin{rem}\label{eqv6} Note that in small dimensions ($m\le 4$) the necessity
condition for homotopy of Volterra q.s.o. is  sufficient as well.
\end{rem}

\begin{cor}\label{eqv66} If $V_0 \sim V_1,$ then for any face
$\Gamma_\alpha$ of the simplex $S^{m-1}$ one has
$V_0\mid_{\Gamma_\alpha} \sim V_1\mid_{\Gamma_\alpha}.$
\end{cor}

\begin{pf} Let $I\setminus\Gamma_\alpha=\{i_1,i_2,\dots,i_k\}$.
According to Theorem \ref{vol1} the restriction of Vorterra q.s.o.
$V$ to $\Gamma_\alpha$, i.e.
$V_\alpha:\Gamma_\alpha\rightarrow\Gamma_\alpha$ is also Volterra
q.s.o., therefore, the corresponding the skew-symmetrical matrix
$A_\alpha$ to $V\mid_{\Gamma_\alpha}$ is a matrix which can be
obtained from the matrix $A_m$ by eliminating the rows
$\{i_1,i_2,\dots,i_m\}$ and the columns $\{i_1,i_2,\dots,i_m\}$.
Now if $V_0 \sim V_1$, then from Theorem \ref{eqv5} it follows
that the operators $V_0,V_1$ lie on the same connected component.
From the definition of the subpffaffinas one concludes that the
operators $V_0\mid_{\Gamma_\alpha}, V_1\mid_{\Gamma_\alpha}$ also
belong to the same connected component, hence again Theorem
\ref{eqv5} implies that $V_0\mid_{\Gamma_\alpha},
V_1\mid_{\Gamma_\alpha}$ are homotopic.
\end{pf}

\begin{cor}\label{eqv7} Let $V_0 \sim V_1$. Then for any face $\Gamma_\alpha$
of the simplex $S^{m-1}$ one has
\begin{itemize}
    \item[(i)] $|Fix(V_0)\cap\Gamma_\alpha|=|Fix(V_1)\cap\Gamma_\alpha|.$
    \item[(ii)] $|Fix(V_0)\cap ri\Gamma_\alpha|=|Fix(V_1)\cap ri\Gamma_\alpha|.$
\end{itemize}
\end{cor}

\begin{pf} (i). Corollary \ref{eqv66} yields that $V_0\mid_{\Gamma_\alpha} \sim
V_1\mid_{\Gamma_\alpha}$ for any $\alpha\subset I$, therefore,
$|Fix(V_0)\cap\Gamma_\alpha|=|Fix(V_0)\cap\Gamma_\alpha|$.

(ii) Now suppose that $\alpha=\{i_1,i_2,\dots,i_k\}$. One can see
that
\begin{equation}\label{dG}
\partial\Gamma_\alpha=\bigcup\limits_{n=1}^k\Gamma_{\alpha_{i_n}},
\end{equation}
here $\alpha_{i_n}=\alpha\setminus\{i_n\}.$ From (i) one finds
that
$|Fix(V_0)\cap\Gamma_{\alpha_{i_n}}|=|Fix(V_0)\cap\Gamma_{\alpha_{i_n}}|$
for any $n=\overline{1,k}$, hence from \eqref{dG} we get
$|Fix(V_0)\cap\partial\Gamma_\alpha|=|Fix(V_0)\cap\partial\Gamma_\alpha|$,
which implies $|Fix(V_0)\cap ri\Gamma_\alpha|=|Fix(V_0)\cap
ri\Gamma_\alpha|.$
\end{pf}

\begin{rem}
The proved corollaries imply that equivalent Volterra q.s.o. have
the same number of fixed points on every face and its interior as
well. Due to this facts one can ask: are there homotopic
extensions of a given Volterra q.s.o. on a face to whole simplex?
Next we are going to study this question.
\end{rem}

Let $\alpha\subset I$, and
$V_0:\Gamma_\alpha\rightarrow\Gamma_\alpha$  be a transitive
Volterra q.s.o. on a face $\Gamma_\alpha$. Denote
\begin{equation}\label{FV}
F_{V_0}(\alpha)=\{V\in {\mathcal{V}}_t^{m-1}:
V\mid_{\Gamma_\alpha}\sim V_0\}.
\end{equation}

\begin{rem}\label{FV1} From the definition of $F_{V_0}$ it follows
for any $V_1,V_2\in F_{V_0}(\a)$ one has  $V_1\mid_{\G_\a}\sim
V_2\mid_{\G_\a}.$
\end{rem}

Since $\Gamma_\a$ is a $(|\alpha|-1)$-dimensional simplex, so
Corollary \ref{eqv66} implies that

\begin{lem}\label{FV1} For any $\alpha,\beta\subset I$ one has
$$\beta\subset\alpha \ \ \ \Rightarrow \ \ \ F_{V_0}(\alpha)\subset F_{V_0}(\beta).$$
In particular, $\forall\alpha\subset I$ one gets
$$F_{V_0}(I)\subset F_{V_0}(\a).$$
\end{lem}

Let
\begin{equation}\label{Vj}
{\mathcal{V}}_t^{m-1}/_\sim=\{{\mathcal{V}}_1,\dots,{\mathcal{V}}_r\}.
\end{equation}

Then from \eqref{FV} and \eqref{Vj} one can see that for any
$V_0\in{\mathcal{V}}_t^{m-1}$ there exits $i\in \{1,2,\dots,r\}$
such that
$$
F_{V_0}(I)={\mathcal{V}}_i.
$$
Therefore, we are interested when $|\alpha|\le n-1$. In this case,
it is clear that any Volterra q.s.o. given on $\G_\a$ can be
extended to a transversal Volterra q.s.o. defined on the simplex
$S^{m-1}$. Note that  such an extension is not unique.

In what follows we shall assume that a Volterra q.s.o. $V_0$ is
defined on the whole simplex $S^{m-1}$, i.e. $V_0\in
{\mathcal{V}}_t^{m-1}$.

\begin{thm}\label{FV2} Let $|\alpha|\le n-1$. Then there are
$i,j\in\{1,2,\dots,r\},$ $i\neq j$ such that
$$
F_{V_0}\cap {\mathcal{V}}_i\neq \emptyset, \ \ \  F_{V_0}\cap
{\mathcal{V}}_j\neq \emptyset
$$
\end{thm}

\begin{pf} Due to $V_0\in{\mathcal{V}}_t^{m-1}/_\sim$ there is  $i\in\{1,2,\dots,r\}$  such that
$V_0\in{\mathcal{V}}_i.$ As before, by $A_m^0$ we denote the
corresponding skew-symmetric matrix. From $|\alpha|\le n-1$ we
have $I\setminus\alpha\neq\emptyset.$ Let $p_0\in
I\setminus\alpha$, $q_0\in I.$  Then $a_{p_0q_0}^0$ is not an
element of $A_\alpha^0= A_m^0\mid_{\G_\a}$. Without loss of
generality we may assume that $a_{p_0q_0}^0>0$. Now we are going
to construct a skew-symmetric matrix
$A_m^{1}=(a_{ij}^1)_{i,j=1}^m$ as follows: if
$i,j\notin\{p_0,q_0\}$, then we put $a_{ij}^1=a_{ij}^0$. We choose
$a_{p_0q_0}^1$ from the segment $[-1,0)$ (Note that
$a_{q_0p_0}^0\in (0,1]$) such that all pffaffians of the matrix
$A_m^1$ is not zero. The existence of such a number comes from
that fact that each pffaffian is a polynomial with respect to
$a_{p_0q_0}^1$ (since all the rest elements are defined),
therefore, its zeros are finite, and such paffaffians are finite
as well. So, all pffaffians are not zero except for finite numbers
of $[-1,0)$. According to the construction $A_m^1$ is a skew
-symmetric, hence the corresponding Volterra q.s.o. $V_1$ is
transversal, and moreover,  $A_m^1\mid_{\G_\a}=A_m^0\mid_{\G\a}$,
i.e. $V_1\mid_{\G_\a}\sim V_0\mid_{\G_\a}$. But $V_1 $ and $V_0$
are not homotopic, since  $a_{p_0q_0}^1$ and $a_{p_0q_0}^0$ have
different sighs, this means that the second order pffaffians have
different sighs too (see Theorem \ref{eqv1}). Let
${\mathcal{V}}_j$ be a set of Volttera operators which are
equivalent to $V_1$. Then the construction shows that $i\neq j$.
\end{pf}

\begin{cor}\label{FV3} Note also that if $|\alpha|\le n-1$, then
\begin{enumerate}
    \item[(i)]  $F_{V_0}(\a)$ is not a subset of any ${\mathcal{V}}_i$
(here $i\in\{1,2,\dots,r\}$);
    \item[(ii)] $F_{V_0}(\a)$ is not a linearly connected set.
\end{enumerate}
\end{cor}

\begin{rem} From the proved Theorem \ref{FV2} we conclude that any
transversal Volterra operator given on a face has not a unique
homotopic extension.
\end{rem}

It is clear that
\begin{equation*}
F_{V_0}(\a)=\bigcup\limits_{i=1}^r\bigl(F_{V_0}(\a)\bigcap{\mathcal{V}}_i\bigr).
\end{equation*}

\begin{thm} If there is some $i\in\{1,2,\dots,r\}$ such that
$F_{V_0}(\a)\bigcap{\mathcal{V}}_i$ is not empty, then it is
linearly connected.
\end{thm}

\begin{pf} Let assume that $F_{V_0}(\a)\bigcap{\mathcal{V}}_i$ is not empty for
some $i\in\{1,2,\dots,r\}$. Then take two elements $V_1,V_2\in
F_{V_0}(\a)\bigcap{\mathcal{V}}_i$.  Now we are going to show such
element can be connected with a path lying in
$F_{V_0}(\a)\bigcap{\mathcal{V}}_i$. Taking into account that
$V_1,V_2\in{\mathcal{V}}_i$ and Theorem \ref{eqv5} we find that
there is a path $\{V_\lambda\}_{\lambda\in [1,2]}\subset
{\mathcal{V}}_i$ connecting them. For any $\lambda\in [1,2]$ one
has $V_\lambda \sim V_1$, hence $V_\lambda\mid_{\G_\a} \sim
V_1\mid_{\G_\a}$. From $V_1\mid_{\G_\a} \sim V_0\mid_{\G_\a}$ we
obtain $V_\lambda\mid_{\G_\a} \sim V_0\mid_{\G_\a}$, this means
that $\{V_\lambda\}_{\lambda\in [1,2]}\subset F_{V_0}(\a)$.
Therefore, $F_{V_0}(\a)\bigcap{\mathcal{V}}_i$ is linearly
connected.\end{pf}

\begin{cor} For any $\alpha\subset I$ one has
\begin{equation*}
F_{V_0}(\a)/_\sim=\bigl\{F_{V_0}(\a)\bigcap{\mathcal{V}}_i\bigr\}_{i=1}^r,
\end{equation*}
here for some $i$ the set $F_{V_0}(\a)\bigcap{\mathcal{V}}_i$ can
be empty. Therefore,
$$|F_{V_0}(\a)/_\sim|\le r.$$
In particular, the equality occurs when $|\alpha|=1$, i.e.
$$F_{V_0}(\a)/_\sim=\bigl\{{\mathcal{V}}_i\bigr\}_{i=1}^r={\mathcal{V}}_t^{m-1}/_\sim.$$
\end{cor}

From Theorem \ref{FV2} we conclude that if two Volterra operators
are homotopic on a face, then there need not be homotopic on the
simplex $S^{m-1}$. But one arises the following problem.

\begin{prob}
How many faces need, on which two Volterra operators are
homotopic, to be homotopic of such operators on the simplex
$S^{m-1}$?
\end{prob}

Now we are going to show that the formulated problem has negative
solution when $m$ is even.

{\bf Example.} Consider a case when $m=4$, then According Theorem
\ref{eqv3} we know that two Volterra q.s.o. are homotopic iff the
signature of corresponding matrices are the same and moreover,
pfaffians of their determinants have the same sign. Let us
consider two transversal Volttera q.s.o. corresponding to the
following matrices
$$A_4^{(1)}=\left(%
\begin{array}{cccc}
  0 & 1 & 1 & 1 \\
  -1 & 0 & 1 & 1 \\
  -1 & -1 & 0 & 1\\
  -1 & -1 & -1 & 0\\
\end{array}%
\right), \ \ \
A_4^{(2)}=\left(%
\begin{array}{cccc}
  0 & \frac12 & 1 & \frac12 \\
  -\frac12 & 0 & \frac12 & 1 \\
  -1 & -\frac12 & 0 & \frac12\\
  -\frac12 & -1 & -\frac12 & 0\\
\end{array}%
\right).$$

Then one can check that the operators $V_1$ and $V_2$ are
homotopic on any proper face of the simplex $S^3$(see Theorem
\ref{eqv2}). Since the pffaffians corresponding to determinants of
the matrices $A_4^{(1)}$ and  $A_4^{(2)}$ are 1 and $-\frac12$,
respectively, therefore due to Theorem \ref{eqv3} we conclude that
they are not homotopic on the whole simplex $S^3$.

\end{document}